\documentclass[10pt,reqno]{amsart}
   \usepackage[centertags]{amsmath}
   \usepackage{amsfonts}
   \usepackage{amsmath}
   \usepackage{amssymb}
   \usepackage{amsthm}
   \usepackage{newlfont}
\usepackage[applemac]{inputenc}

\usepackage{multirow, bigdelim}

\usepackage{epsfig}
\usepackage{pst-all}
\usepackage{pstricks}

\usepackage{pst-pdgr}%
\usepackage{caption}
\usepackage{subcaption}

\usepackage{color}

\usepackage{graphicx}

\usepackage{bm}

\allowdisplaybreaks[3]

\theoremstyle{plain}
\newtheorem{theorem}{Theorem}[section]

\theoremstyle{definition}

\theoremstyle{remark}
\newtheorem{remark}[theorem]{Remark}
\newtheorem*{remark*}{Remark}

\numberwithin{equation}{section}

\newcommand\QQ{{\mathbb Q}}

\title[QBD and Jacobi-Koornwinder bivariate polynomials]{QBD processes associated with Jacobi-Koornwinder \\ bivariate polynomials and urn models}

\author{Lidia Fern\'andez and Manuel D. de la Iglesia}

\address{Lidia Fern\'andez\\
IEMATH-GR and Departamento de Matem\'atica Aplicada\\
Universidad de Granada\\
18071, Granada, Spain.}
\email{lidiafr@ugr.es}
\address{Manuel D. de la Iglesia\\
Instituto de Matem\'aticas \\
Universidad Nacional Aut\'onoma de M\'exico \\
Circuito Exterior, C.U.\\
04510, Mexico D.F. Mexico.}
\email{mdi29@im.unam.mx}
\thanks{The work of the first author was partially supported by FEDER/Junta de Andaluc\'ia under the research project A-FQM-246-UGR20; MCIN/AEI 10.13039/501100011033 and FEDER funds by PGC2018-094932-B-I00; and IMAG-Mar\'ia de Maeztu grant CEX2020-001105-M.
The work of the second author was partially supported by PAPIIT-DGAPA-UNAM grant IN106822 (M\'exico) and CONACYT grant A1-S-16202 (M\'exico).}

\date{\today}



\begin{document}

\maketitle

\begin{abstract}
We study a family of quasi-birth-and-death (QBD) processes associated with the so-called first family of Jacobi-Koornwinder bivariate polynomials. These polynomials are orthogonal on a bounded region typically known as the swallow tail. We will explicitly compute the coefficients of the three-term recurrence relations generated by these QBD polynomials and study the conditions under we can produce families of discrete-time QBD processes. Finally, we show an urn model associated with one special case of these QBD processes.
\end{abstract}

\section{Introduction}

In last few years there has been an increasing activity in the study of the spectral representation of \emph{quasi-birth-and-death} (QBD) \emph{processes}, extending the pioneering work of S. Karlin and J. McGregor \cite{KMc2, KMc3, KMc6} in the 1950s (see also the recent monograph \cite{MDIB}). These processes are a natural extension of the so-called birth-death chains, where the state space, instead of $\mathbb{N}_0$, is given by pairs of the form $(n,k)$, where $n\in\mathbb{N}_0$ is usually called the \emph{level}, while $1\leq k\leq r_n$ is referred to as the \emph{phase} (which may depend on the different levels). For a general setup see \cite{LaR}. The transition probability matrix (discrete-time) or the infinitesimal operator matrix (continuous-time) of the QBD process is then block tridiagonal. If $r_n=1$ for all $n\in\mathbb{N}_0$ then we go back to classical birth-death chains. If $r_n=N$ for all $n\in\mathbb{N}_0$, where $N$ is a positive integer, then all blocks in the Jacobi matrix have the same dimension $N\times N$. In this case, the spectral analysis can be performed by using \emph{matrix-valued orthogonal polynomials} (see \cite{DRSZ, G2} for the discrete-time case and \cite{DR} for the continuous-time case). Many examples have been analyzed in this direction by using spectral methods in the last few years (see \cite{ Clay, DRSZ, G2, GdI3, dI1, dIJ, dIR}).

A natural source of examples of more complicated QBD processes comes from the theory of \emph{multivariate orthogonal polynomials} (of dimension $d$), where now the number of phases is given by $r_n=\binom{n+d-1}{n}$. In \cite{FdI} we performed the spectral analysis in the general setting of this situation as well as obtained results about recurrence and the invariant measure of these processes in terms of the spectral measure supported on some region $\Omega\subset\mathbb{R}^d$. We also applied our results to several examples of bivariate orthogonal polynomials ($d=2$), namely product orthogonal polynomials, orthogonal polynomials on a parabolic domain and orthogonal polynomials on the triangle. The aim of this paper is to continue our previous work but now we will focus on the so-called first family of bivariate Jacobi-Koornwinder polynomials (see \cite[Section 2.7]{DX14}), first introduced by T. Koornwinder in \cite{Ko74} (see also the review paper \cite{Ko75} where they are called Class VI). These polynomials are supported in the so-called swallow tail region (see Figure \ref{fig1}) and they are eigenfunctions of two independent differential operators of orders two and four. Some properties such as a Rodrigues-type expression or an expansion in terms of James-type zonal polynomials can be found in \cite{Sp76} and \cite{KS78}. They are considered a highly non-trivial generalization of the Jacobi polynomials. Yuan Xu proved some cubature rules for specific values of the parameters \cite{Xu12,Xu17}.

The paper is organized as follows. In Section \ref{sec2} we introduce the Jacobi-Koornwinder polynomials we will be working with. Then we normalize the polynomials in such a way that they are equal to 1 at one of the corners of the swallow tail region (specifically at the point $(1,1)$). With this family of polynomials we derive the coefficients of the two three-term recurrence relations (one for each variable) in terms of the coefficients of the three-term recurrence relation of the classical Jacobi polynomials on $[0,1]$. In Section \ref{sec3} we will study under what conditions we may provide a probabilistic interpretation of the linear convex combination of the two Jacobi matrices associated with the three-term recurrence relations. Under these conditions we compute the Karlin-McGregor formula, the invariant measure and study recurrence of the family of discrete-time QBD processes. Finally, in Section \ref{sec4}, we give an urn model associated with one of the QBD processes introduced in Section \ref{sec3}, for the special case of $\beta=\alpha$.

\section{Bivariate Jacobi-Koornwinder polynomials}\label{sec2}

\begin{figure}[h]
\includegraphics[scale=0.5]{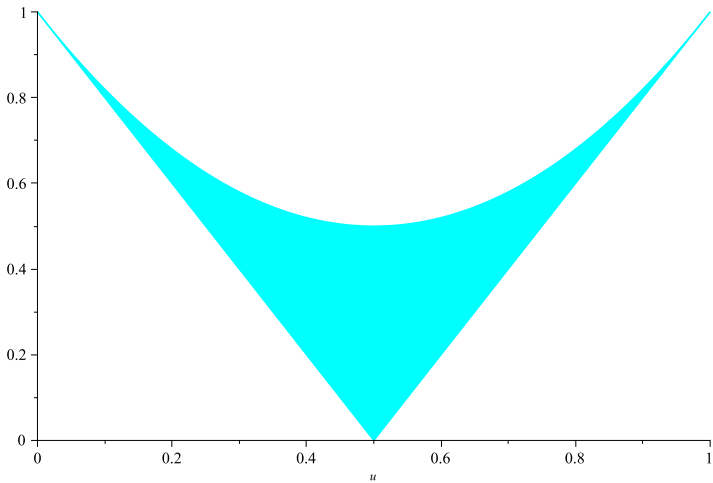}
\caption{Swallow tail region $\Omega$ where the weight function is defined.}
\label{fig1}
\end{figure}

In \cite[Section 2.7]{DX14} and \cite{Ko75}, the Jacobi-Koornwinder polynomials are constructed in terms of the Jacobi weight function supported on $[-1,1]$. The swallow tail region $\Omega$ is then contained in the bounded rectangle $[-2,2]\times[-1,1]$. For convenience we consider a change of variables ($u\mapsto 4u-2, v\mapsto 2v-1$) such that the swallow tail region $\Omega$ is contained inside the unit square $[0,1]\times[0,1]$. Then the region $\Omega$ is given by (see Figure \ref{fig1})
$$
\Omega=\left\{(u,v)\,:\,2u+v-1>0, 1-2u+v>0, 2u^2-2u-v+1>0 \right\}.
$$

The weight function acting on this region will be given by
\begin{equation}\label{wKpJ2}
W_{\alpha,\beta,\gamma}(u,v)=\frac{1}{C}(1-2u+v)^{\alpha} (2u+v-1)^{\beta} (2u^2-2u-v+1)^{\gamma},
\end{equation}
where $C$ is the normalizing constant
$$
C=\frac{2^{\alpha+\beta-\gamma+2}\Gamma(\alpha+1)\Gamma(\beta+1)\Gamma(\gamma+1)\Gamma(2\alpha+2\gamma+2)\Gamma(2\beta+2\gamma+2)\Gamma(\alpha+\beta+\gamma+3)}{\Gamma(\alpha+\gamma+1)\Gamma(\beta+\gamma+1)\Gamma(\alpha+\beta+2\gamma+3)\Gamma(2\alpha+2\beta+2\gamma+5)},
$$
such that $\int_{\Omega}W_{\alpha,\beta,\gamma}(u,v)dudv=1$. To ensure integrability we need to have $\alpha,\beta,\gamma>-1$,  $\alpha+\gamma+3/2>0$ and $\beta+\gamma+3/2>0$. This normalized constant was computed in \cite[Lema 6.1]{Sp76}.

As it was pointed out in \cite[Proposition 2.7.3]{DX14} the monic Jacobi-Koornwinder polynomials $P_{n,k}^{\alpha,\beta,\gamma}(u,v)$ satisfy the following second-order partial differential equation (after the change of variables):
\begin{equation*}\label{pdeq}
\mathcal{D}_{\alpha,\beta,\gamma}P_{n,k}^{\alpha,\beta,\gamma}(u,v)=-\lambda_{n,k}^{\alpha,\beta,\gamma}P_{n,k}^{\alpha,\beta,\gamma}(u,v),
\end{equation*}
where
\begin{align*}
\mathcal{D}_{\alpha,\beta,\gamma}&=\left[u(1-u)-(1-v)/4\right]\partial_{uu}+(1-v)(2u-1)\partial_{uv}+\left[(2u-1)^2+v(1-2v)\right]\partial_{vv}\\
&\quad+\left[\beta+\gamma+3/2-(\alpha+\beta+2\gamma+3)u\right]\partial_u+2\left[(\beta-\alpha)u-(\alpha+\beta+\gamma+5/2)v+\alpha+1\right]\partial_v,\\
\lambda_{n,k}^{\alpha,\beta,\gamma}&=n(n+\alpha+\beta+2\gamma+2)+k(k+\alpha+\beta+1).
\end{align*}
With this partial differential equation it is possible to generate all the monic Jacobi-Koornwinder polynomials for any values of $\alpha, \beta, \gamma$. For the special cases of $\gamma=\pm1/2$ it is possible to write them in terms of classical Jacobi polynomials (see \cite[Proposition 2.7.2]{DX14}). Another way to compute the Jacobi-Koornwinder polynomials is by using the Rodrigues-type formula found in \cite[Section 5]{Sp76}.

\medskip

Let us now introduce a new set of polynomials $Q_{n,k}^{\alpha,\beta,\gamma}(u,v)$ normalized in such a way that $Q_{n,k}^{\alpha,\beta,\gamma}(1,1)=1$ for all $n\geq0$ and $0\leq k\leq n$. $Q_{n,k}^{\alpha,\beta,\gamma}(u,v)$ can also be defined as
$$
Q_{n,k}^{\alpha,\beta,\gamma}(u,v)=\sigma_{n,k}^{-1}P_{n,k}^{\alpha,\beta,\gamma}(u,v),
$$
where $\sigma_{n,k}=P_{n,k}^{\alpha,\beta,\gamma}(1,1)$ is given by
\begin{equation}\label{sigmm}
\sigma_{n,k}=\frac{2^{2k-n+1}(2\gamma+2)_{n-k-1}(\alpha+1)_k(\alpha+\gamma+3/2)_n}{(\gamma+3/2)_{n-k-1}(k+\alpha+\beta+1)_k(n+k+\alpha+\beta+2\gamma+2)_{n-k}(n+\alpha+\beta+\gamma+3/2)_k}.
\end{equation}
Here we are using the standard notation for the Pochhammer symbol $(a)_0=1, (a)_n=a(a+1)\cdots(a+n-1), n\geq1$. The vector polynomials $\QQ_n=(Q_{n,0}^{\alpha,\beta,\gamma},Q_{n,1}^{\alpha,\beta,\gamma},\dots, Q_{n,n}^{\alpha,\beta,\gamma}), n\geq0,$ satisfy the three-term recurrence relations
\begin{equation}\label{TTRRQ}
\begin{aligned}
u\, {\mathbb Q}_n(u,v) & = A_{n,1} {\mathbb Q}_{n+1}(u,v)+ B_{n,1} {\mathbb Q}_n(u,v) + C_{n,1} {\mathbb Q}_{n-1}(u,v), \\
v\, {\mathbb Q}_n(u,v) & = A_{n,2} {\mathbb Q}_{n+1}(u,v)+ B_{n,2} {\mathbb Q}_n(u,v) + C_{n,2} {\mathbb Q}_{n-1}(u,v).
\end{aligned}
\end{equation}
It is possible to compute explicitly the coefficients $A_{n,i}, B_{n,i},C_{n,i}, i=1,2,$ using \cite[Section 9]{Sp76} and \eqref{sigmm}.  For that let us introduce the following notation:
\begin{equation}\label{coefTTRR}
\begin{split}
a_n&=\frac{(n+\alpha+1)(n+\alpha+\beta+1)}{(2n+\alpha+\beta+1)(2n+\alpha+\beta+2)},\\
b_n&=\frac{(n+\alpha+1)(n+1)}{(2n+\alpha+\beta+1)(2n+\alpha+\beta+2)}+\frac{(n+\beta)(n+\alpha+\beta)}{(2n+\alpha+\beta)(2n+\alpha+\beta+1)},\\
c_n&=\frac{n(n+\beta)}{(2n+\alpha+\beta)(2n+\alpha+\beta+1)},\\
\delta_{n,k}&=(n-k)(n+k+\alpha+\beta+1).
\end{split}
\end{equation}
Observe that $a_n,b_n,c_n$ are the coefficients of the three-term recurrence relation satisfied by the classical Jacobi polynomials $Q_n^{(\beta,\alpha)}(x)$ on $[0,1]$ normalized by $Q_n^{(\beta,\alpha)}(1)=1$ (see \cite[Section 5]{GdI4} for instance). In particular we always have that $a_n,c_n>0, b_n\geq0,$ and $a_n+b_n+c_n=1$, i.e. they are probabilities. The norms of these Jacobi polynomials (which will be used later) are given by
\begin{equation}\label{normJ}
\|Q_n^{(\beta,\alpha)}\|_w^2=\frac{\Gamma(\alpha+1)\Gamma(\alpha+\beta+2)\Gamma(n+1)\Gamma(n+\beta+1)}{\Gamma(\beta+1)\Gamma(n+\alpha+1)\Gamma(n+\alpha+\beta+1)(2n+\alpha+\beta+1)},
\end{equation}
where $w$ is the normalized Jacobi weight (see (5.2) of \cite{GdI4}).

On one side, the matrices $A_{n,1}$, $B_{n,1}$ and $C_{n,1}$ in \eqref{TTRRQ} are of the form
\begin{equation}\label{ABC1}
\begin{array}{c}
A_{n,1}=\left[
\begin{array}{ccccccc}
a_{n,0} &            &        &             & 0  \\
        & a_{n,1}    &        &             & \vdots \\
        &            & \ddots &             &         \\
        &            &        & a_{n,n}   & 0
\end{array}
\right], \quad
C_{n,1}=\left[
\begin{array}{cccccccc}
c_{n,0}     &            &        &               \\
            & c_{n,1}    &        &               \\
            &            & \ddots &             \\
            &            &        & c_{n,n-1} \\
     0      &            & \dots   &  0
\end{array}
\right],
\\[1cm]
B_{n,1}=\left[
\begin{array}{ccccccc}
b_{n,0} & e_{n,0}&         &        &  \\
d_{n,1} & b_{n,1}& e_{n,1}&        &           \\
              & \ddots & \ddots  & \ddots &           \\
              &        & d_{n,n-1} &  b_{n,n-1}& e_{n,n-1} \\
              &        &       &  d_{n,n}     & b_{n,n}
\end{array}
\right],
\end{array}
\end{equation}
where the entries of $A_{n,1}$, $B_{n,1}$ and $C_{n,1}$ (see \eqref{coefTTRR}) are given by
\begin{equation}\label{coeffpar0}
\begin{split}
a_{n,k}&=\frac{1}{2}a_{n+\gamma+\frac{1}{2}}\frac{\delta_{n+2\gamma+1,k}}{\delta_{n+\gamma+\frac{1}{2},k}},\quad k=0,1,\ldots,n, \\
c_{n,k}&=\frac{1}{2}c_{n+\gamma+\frac{1}{2}}\frac{\delta_{n,k}}{\delta_{n+\gamma+\frac{1}{2},k}},\quad k=0,1,\ldots,n-1,
\\
e_{n,k}&=\frac{1}{2}a_{k}\frac{\delta_{n+\gamma+\frac{1}{2},k+\gamma+\frac{1}{2}}}{\delta_{n+\gamma+\frac{1}{2},k}},\quad k=0,1,\ldots,n-1,\\
d_{n,k}&=\frac{1}{2}c_{k}\frac{\delta_{n+\gamma+\frac{1}{2},k-\gamma-\frac{1}{2}}}{\delta_{n+\gamma+\frac{1}{2},k}},\quad k=1,2,\ldots,n-1,\\
b_{n,k}&=1-a_{n,k}-c_{n,k}-d_{n,k}-e_{n,k},\quad k=0,1,\ldots,n.
\end{split}
\end{equation}
\begin{remark}
The coefficient $b_{n,k}$ can also be written as
$$
b_{n,k}=\frac{1}{2}(b_{n+\gamma+1/2}+b_k)+\frac{1-4\gamma^2}{4(\beta^2-\alpha^2)}(2b_{n+\gamma+1/2}-1)(2b_{k}-1).
$$
Also, from \eqref{coefTTRR} it is possible to see that
$$
\frac{\delta_{n+2\gamma+1,k}}{\delta_{n+\gamma+\frac{1}{2},k}}+\frac{\delta_{n,k}}{\delta_{n+\gamma+\frac{1}{2},k}}+\frac{\delta_{n+\gamma+\frac{1}{2},k-\gamma-\frac{1}{2}}}{\delta_{n+\gamma+\frac{1}{2},k}}+\frac{\delta_{n+\gamma+\frac{1}{2},k+\gamma+\frac{1}{2}}}{\delta_{n+\gamma+\frac{1}{2},k}}=4.
$$
\end{remark}
On the other side, the matrices $A_{n,2}$, $B_{n,2}$ and $C_{n,2}$ are tridiagonal matrices of the form
\begin{equation}\label{ABC2}
\begin{array}{c}
A_{n,2}=\left[
\begin{array}{ccccccc}
a_{n,0}^{(2)} & a_{n,0}^{(3)}&         &          &  \\
a_{n,1}^{(1)} & a_{n,1}^{(2)}& a_{n,1}^{(3)} &         &    \\
        & \ddots & \ddots  & \ddots  &         \\
 &        & a_{n,n}^{(1)} & a_{n,n}^{(2)} & a_{n,n}^{(3)}
\end{array}
\right],
\quad
B_{n,2}=\left[
\begin{array}{ccccccc}
b_{n,0}^{(2)} & b_{n,0}^{(3)}&         &        &  \\
b_{n,1}^{(1)} & b_{n,1}^{(2)}& b_{n,1}^{(3)} &        &           \\
              & \ddots & \ddots  & \ddots &           \\
              &        & b_{n,n-1}^{(1)} &  b_{n,n-1}^{(2)} & b_{n,n-1}^{(3)}  \\
              &        &       &  b_{n,n}^{(1)}     & b_{n,n}^{(2)}
\end{array}
\right],
\\[1cm]
C_{n,2}=\left[
\begin{array}{ccccccc}
c_{n,0}^{(2)} & c_{n,0}^{(3)}&         &          &  \\
c_{n,1}^{(1)} & c_{n,1}^{(2)}& c_{n,1}^{(3)} &         &    \\
        & \ddots & \ddots  & \ddots  &         \\
         &        & c_{n,n-2}^{(1)} & c_{n,n-2}^{(2)} & c_{n,n-2}^{(3)} \\
         &        &        &  c_{n,n-1}^{(1)} & c_{n,n-1}^{(2)} \\
        &        &        &         & c_{n,n}^{(1)}
\end{array}
\right],
\end{array}
\end{equation}
where the entries of $A_{n,2}$, $B_{n,2}$ and $C_{n,2}$ (see again \eqref{coefTTRR}) are given by
\begin{equation}\label{coeffpar}
\begin{split}
a_{n,k}^{(1)}&=2c_ka_{n+\gamma+\frac{1}{2}}\frac{\delta_{n+\gamma+\frac{1}{2},k-\gamma-\frac{1}{2}}\delta_{n+\gamma+1,k-\gamma-1}}{\delta_{n+\frac{1}{2}(\gamma+\frac{1}{2}),k-\frac{1}{2}(\gamma+\frac{1}{2})}\delta_{n+\frac{1}{2}(\gamma+\frac{3}{2}),k-\frac{1}{2}(\gamma+\frac{3}{2})}},\quad k=1,\ldots,n,\\
a_{n,k}^{(2)}&=(2b_k-1)a_{n+\gamma+\frac{1}{2}}\frac{\delta_{n+2\gamma+1,k}}{\delta_{n+\gamma+\frac{1}{2},k}},\quad  k=0,1,\ldots,n,\\
a_{n,k}^{(3)}&=2a_ka_{n+\gamma+\frac{1}{2}}\frac{\delta_{n+\gamma+\frac{1}{2},k+\gamma+\frac{1}{2}}\delta_{n+\gamma+1,k+\gamma+1}}{\delta_{n+\frac{1}{2}(\gamma+\frac{1}{2}),k+\frac{1}{2}(\gamma+\frac{1}{2})}\delta_{n+\frac{1}{2}(\gamma+\frac{3}{2}),k+\frac{1}{2}(\gamma+\frac{3}{2})}},\quad k=0,1,\ldots,n,
\\
b_{n,k}^{(1)}&=(2b_{n+\gamma+\frac{1}{2}}-1)c_k\frac{\delta_{n+\gamma+\frac{1}{2},k-\gamma-\frac{1}{2}}}{\delta_{n+\gamma+\frac{1}{2},k}},\quad k=1,\ldots,n,
\\
b_{n,k}^{(2)}&=1-b_{n,k}^{(1)}-b_{n,k}^{(3)}-a_{n,k}^{(1)}-a_{n,k}^{(2)}-a_{n,k}^{(3)}-c_{n,k}^{(1)}-c_{n,k}^{(2)}-c_{n,k}^{(3)},\quad k=0,1,\ldots,n,
\\
b_{n,k}^{(3)}&=(2b_{n+\gamma+\frac{1}{2}}-1)a_k\frac{\delta_{n+\gamma+\frac{1}{2},k+\gamma+\frac{1}{2}}}{\delta_{n+\gamma+\frac{1}{2},k}},\quad k=0,1,\ldots,n-1,
\\
c_{n,k}^{(1)}&=2c_kc_{n+\gamma+\frac{1}{2}}\frac{\delta_{n,k}\delta_{n-\frac{1}{2},k-\frac{1}{2}}}{\delta_{n+\frac{1}{2}(\gamma-\frac{1}{2}),k+\frac{1}{2}(\gamma-\frac{1}{2})}\delta_{n+\frac{1}{2}(\gamma+\frac{1}{2}),k+\frac{1}{2}(\gamma+\frac{1}{2})}},\quad k=1,\ldots,n,
\\
c_{n,k}^{(2)}&=(2b_k-1)c_{n+\gamma+\frac{1}{2}}\frac{\delta_{n,k}}{\delta_{n+\gamma+\frac{1}{2},k}},\quad k=0,1,\ldots,n-1,\\
c_{n,k}^{(3)}&=2a_kc_{n+\gamma+\frac{1}{2}}\frac{\delta_{n,k}\delta_{n-\frac{1}{2},k+\frac{1}{2}}}{\delta_{n+\frac{1}{2}(\gamma-\frac{1}{2}),k-\frac{1}{2}(\gamma-\frac{1}{2})}\delta_{n+\frac{1}{2}(\gamma+\frac{1}{2}),k-\frac{1}{2}(\gamma+\frac{1}{2})}},\quad  k=0,1,\ldots,n-2.
\end{split}
\end{equation}

The coefficients for $A_{n,i}, B_{n,i},C_{n,i}, i=1,2,$ can be significantly simplified for the values of $\gamma=\pm1/2$. For $\gamma=-1/2$ we get
$$
a_{n,k}=\frac{1}{2}a_n,\quad c_{n,k}=\frac{1}{2}c_n,\quad e_{n,k}=\frac{1}{2}a_k,\quad d_{n,k}=\frac{1}{2}c_k,\quad b_{n,k}=\frac{1}{2}(b_n+b_k),
$$
and
\begin{align*}
a_{n,k}^{(1)}&=2a_nc_k,\quad a_{n,k}^{(2)}=a_n(2b_k-1),\quad a_{n,k}^{(3)}=2a_na_k,\\
b_{n,k}^{(1)}&=(2b_n-1)c_k,\quad b_{n,k}^{(2)}=\frac{1}{2}\left(1+(2b_n-1)(2b_k-1)\right),\quad b_{n,k}^{(3)}=(2b_n-1)a_k,\\
c_{n,k}^{(1)}&=2c_nc_k,\quad c_{n,k}^{(2)}=c_n(2b_k-1),\quad c_{n,k}^{(3)}=2c_na_k,
\end{align*}
while for $\gamma=1/2$ we obtain
\begin{align*}
a_{n,k}&=\frac{1}{2}a_{n+1}\frac{\delta_{n+2,k}}{\delta_{n+1,k}},\quad c_{n,k}=\frac{1}{2}c_{n+1}\frac{\delta_{n,k}}{\delta_{n+1,k}},\quad e_{n,k}=\frac{1}{2}a_k\frac{\delta_{n+1,k+1}}{\delta_{n+1,k}},\\
d_{n,k}&=\frac{1}{2}c_k\frac{\delta_{n+1,k-1}}{\delta_{n+1,k}},\quad b_{n,k}=\frac{1}{2}(b_{n+1}+b_k),
\end{align*}
and
\begin{align*}
a_{n,k}^{(1)}&=2a_{n+1}c_k\frac{\delta_{n+\frac{3}{2},k-\frac{3}{2}}}{\delta_{n+\frac{1}{2},k-\frac{1}{2}}},\quad a_{n,k}^{(2)}=a_{n+1}(2b_k-1)\frac{\delta_{n+2,k}}{\delta_{n+1,k}},\quad a_{n,k}^{(3)}=2a_{n+1}a_k\frac{\delta_{n+\frac{3}{2},k+\frac{3}{2}}}{\delta_{n+\frac{1}{2},k+\frac{1}{2}}},\\
b_{n,k}^{(1)}&=(2b_{n+1}-1)c_k\frac{\delta_{n+1,k-1}}{\delta_{n+1,k}},\quad b_{n,k}^{(2)}=\frac{1}{2}\left(1+(2b_{n+1}-1)(2b_k-1)\right),\quad b_{n,k}^{(3)}=(2b_{n+1}-1)a_k\frac{\delta_{n+1,k+1}}{\delta_{n+1,k}},\\
c_{n,k}^{(1)}&=2c_{n+1}c_k\frac{\delta_{n-\frac{1}{2},k-\frac{1}{2}}}{\delta_{n+\frac{1}{2},k+\frac{1}{2}}},\quad c_{n,k}^{(2)}=c_{n+1}(2b_k-1)\frac{\delta_{n,k}}{\delta_{n+1,k}},\quad c_{n,k}^{(3)}=2c_{n+1}a_k\frac{\delta_{n-\frac{1}{2},k+\frac{1}{2}}}{\delta_{n+\frac{1}{2},k-\frac{1}{2}}}.
\end{align*}

\begin{remark}\label{rem32}
The normalization of the polynomials $Q_{n,k}^{\alpha,\beta,\gamma}(u,v)$ such that $Q_{n,k}^{\alpha,\beta,\gamma}(1,1)=1$
will guarantee us that the sum of all rows of the corresponding Jacobi matrices $J_1$ and $J_2$ (see \eqref{JacMat} below) is exactly 1. This does not mean that both $J_1$ and $J_2$ are stochastic matrices or have some probabilistic interpretation, something that we will discuss in the next section. We could have used another ``corner'' of the region $\Omega$ (see Figure \ref{fig1}) like $(0,1)$ or $(1/2,0)$. On one side, it turns out that normalization at the point $(1/2,0)$ will not provide us Jacobi matrices with probabilistic interpretation. On the other side, normalization at the point $(0,1)$ is somehow ``symmetric'' to the normalization at the point $(1,1)$. Indeed, we have that $P_{n,k}^{\alpha,\beta,\gamma}(0,1)=(-1)^{n+k}\sigma_{n,k}$, where $\sigma_{n,k}$ is given by \eqref{sigmm}, and the corresponding new vector polynomials $\widetilde\QQ_n$ satisfy the three-term recurrence relations
\begin{equation*}
\begin{aligned}
-u\, \widetilde{\mathbb Q}_n(u,v) & = \widetilde A_{n,1} \widetilde{\mathbb Q}_{n+1}(u,v)+ \widetilde B_{n,1} \widetilde{\mathbb Q}_n(u,v) + \widetilde C_{n,1} \widetilde{\mathbb Q}_{n-1}(u,v), \\
v\,\widetilde{\mathbb Q}_n(u,v) & = \widetilde A_{n,2} \widetilde{\mathbb Q}_{n+1}(u,v)+ \widetilde B_{n,2} \widetilde{\mathbb Q}_n(u,v) + \widetilde C_{n,2} \widetilde{\mathbb Q}_{n-1}(u,v),
\end{aligned}
\end{equation*}
where the coefficients $\widetilde A_{n,i}, \widetilde B_{n,i}, \widetilde C_{n,i}, i=1,2,$ are \emph{exactly the same} as the coefficients $A_{n,i}, B_{n,i},C_{n,i}, i=1,2,$ but changing $\alpha$ by $\beta$ and $\beta$ by $\alpha$, except for $\widetilde B_{n,1}$ where we have $\widetilde B_{n,1}=B_{n,1}-I$ (changing $\alpha$ by $\beta$ again and viceversa). In this case we have that the sum of the rows of the Jacobi matrix $J_1$ is 0, while the sum of the rows of the Jacobi matrix $J_2$ is 1. For more comments about the choice of normalizing corners the reader can consult \cite[Section 6]{FdI}.

\end{remark}

\section{QBD processes associated with Jacobi-Koornwinder bivariate polynomials}\label{sec3}

\begin{figure}[h]
\includegraphics[scale=0.65]{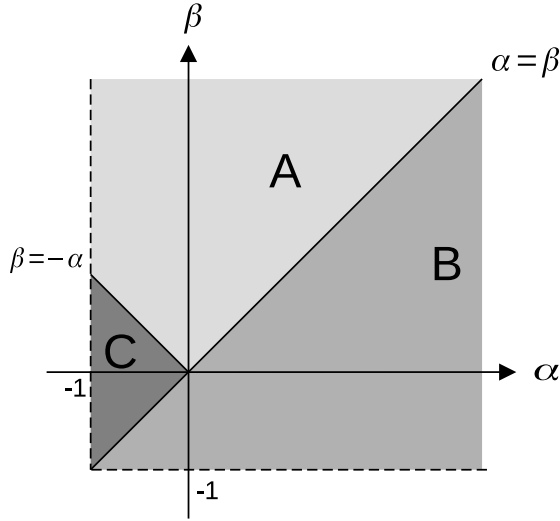}
\caption{The regions \textbf{A}, \textbf{B} and \textbf{C} where $\tau$ may take different values in order to have a stochastic matrix $\bm P$ (courtesy of C. Juarez).}
\label{fig3}
\end{figure}

In this section we will study under what conditions we may provide a probabilistic interpretation of the coefficients of the three-term recurrence relations \eqref{ABC1} and \eqref{ABC2}. From the recurrence relations \eqref{TTRRQ} we can define the following two block tridiagonal Jacobi matrices
\begin{equation}\label{JacMat}
J_{1}=\left(\begin{array}{cccccc}
B_{0,1}   & A_{0,1}   &         &       &\bigcirc \\
C_{1,1}   &B_{1,1}   & A_{1,1} &       &          \\
          &C_{2,1}   & B_{2,1} &A_{2,1} &         \\
  \bigcirc         &           & \ddots  & \ddots & \ddots
\end{array}\right),\quad J_{2}=\left(\begin{array}{cccccc}
B_{0,2}   & A_{0,2}   &         &       &\bigcirc \\
C_{1,2}   &B_{1,2}   & A_{1,2} &       &          \\
          &C_{2,2}   & B_{2,2} &A_{2,2} &         \\
  \bigcirc         &           & \ddots  & \ddots & \ddots
\end{array}\right).
\end{equation}
By construction of the vector polynomials $\mathbb{Q}_n$ (see \eqref{TTRRQ}) we always have that $J_i\bm e=\bm e, i=1,2$, where $\bm e=\left(1, 1, 1, \ldots\right)^T$ is the semi-infinite vector with all components equal to 1. We now consider the linear convex combination of $J_1$ and $J_2$ in the following way
\begin{equation}\label{Pst}
\bm P=(1-\tau)J_1+\tau J_2,\quad 0\leq\tau\leq1.
\end{equation}
We would like to see under what conditions we get a probabilistic interpretation of $\bm P$. In particular we will see when $\bm P$ is a stochastic matrix. We immediately have that $\bm P\bm e=\bm e$ but now we need all entries of $\bm P$ to be positive (except possibly for the main block diagonal, where we only need to be nonnegative). Therefore, looking at the nonzero entries of $\bm P$, we need to have
\begin{align*}
\tau a_{n,k}^{(1)}&>0,\quad k=1,\ldots,n,\quad n\geq1,\\
(1-\tau)a_{n,k}+\tau a_{n,k}^{(2)}&>0,\quad k=0,1,\ldots,n,\quad n\geq0,\\
\tau a_{n,k}^{(3)}&>0,\quad k=0,1,\ldots,n,\quad n\geq0,\\
(1-\tau)d_{n,k}+\tau b_{n,k}^{(1)}&\geq0,\quad k=1,\ldots,n,\quad n\geq1,\\
(1-\tau)b_{n,k}+\tau b_{n,k}^{(2)}&\geq0,\quad k=0,1,\ldots,n,\quad n\geq0,\\
(1-\tau)e_{n,k}+\tau b_{n,k}^{(3)}&\geq0,\quad k=0,1,\ldots,n-1,\quad n\geq1,\\
\tau c_{n,k}^{(1)}&>0,\quad k=1,\ldots,n,\quad n\geq1,\\
(1-\tau)c_{n,k}+\tau c_{n,k}^{(2)}&>0,\quad k=1,\ldots,n,\quad n\geq1,\\
\tau c_{n,k}^{(3)}&>0,\quad k=0,1,\ldots,n-2,\quad n\geq2.
\end{align*}
From the definition of the coefficients \eqref{coeffpar0} and \eqref{coeffpar} we have the following properties
\begin{align*}
a_{n,k}-a_{n,k}^{(2)}&=a_{n,k}(3-4b_k),\quad c_{n,k}-c_{n,k}^{(2)}=c_{n,k}(3-4b_k),\\
d_{n,k}-b_{n,k}^{(1)}&=d_{n,k}(3-4b_{n+\gamma+1/2}),\quad e_{n,k}-b_{n,k}^{(3)}=e_{n,k}(3-4b_{n+\gamma+1/2}),
\end{align*}
where $b_n$ is defined by \eqref{coefTTRR}. First, from \eqref{coeffpar0}, we have that $a_{n,k}, c_{n,k}>0$ and $ d_{n,k}, e_{n,k}\geq0$ as long as $\gamma>-1$ for $\alpha,\beta\geq-1/2$ and $\gamma+3/2\geq\max\{-\alpha,-\beta\}$ in any other case. Under these conditions on $\gamma$ we have, for instance for the case $(1-\tau)a_{n,k}+\tau a_{n,k}^{(2)}>0$, that the parameter $\tau$ must be chosen so that
$$
\tau<\frac{1}{3-4b_k},\quad k=0,1,\ldots n.
$$
Similar considerations can be made for the inequalities $(1-\tau)c_{n,k}+\tau c_{n,k}^{(2)}>0$, $(1-\tau)d_{n,k}+\tau b_{n,k}^{(1)}>0$ and $(1-\tau)e_{n,k}+\tau b_{n,k}^{(3)}>0$. That means that the behavior of $(3-4b_n)^{-1}$ as $n\geq0$ will be the main ingredient in order to find these upper bounds for the parameter $\tau$. In particular, we have
\begin{equation}\label{partCn}
\frac{\partial}{\partial n}(3-4b_n)^{-1}=\frac{8(\alpha^2- \beta^2)(2n+\alpha+\beta+1)}{(3\alpha^2+ 2\alpha\beta+4\alpha n-\beta^2+4\beta n+4n^2+2\alpha+2\beta+4n)^2},
\end{equation}
which behavior, for different values of $\alpha$ and $\beta,$ will be the key to analyze the inequalities given above. Therefore let us define the following $\gamma$-dependent constant
\begin{equation*}\label{cgam}
C_\gamma=\frac{1}{3-4b_{\gamma+\frac{1}{2}}}=\frac{(\alpha+\beta+2\gamma+1)(\alpha+\beta+2\gamma+3)}{(\alpha+\beta+2\gamma+1)(\alpha+\beta+2\gamma+3)+2(\alpha^2-\beta^2)}.
\end{equation*}
For $\gamma=\pm1/2$ we have
$$
C_{-1/2}=\frac{\alpha+\beta+2}{3\alpha-\beta+2},\quad C_{1/2}=\frac{(\alpha+\beta+2)(\alpha+\beta+4)}{(\alpha+\beta+2)(\alpha+\beta+4)+2(\alpha^2-\beta^2)}.
$$
Initially we have that $\alpha,\beta,\gamma>-1$. We will divide the two-dimensional region $\alpha,\beta>-1$ in 3 parts, \textbf{A}, \textbf{B} and \textbf{C}, and then study $\gamma$ for each of these regions (see Figure \ref{fig3}). After some extensive computations using \eqref{partCn} we have the following:

\begin{itemize}
\item In the region \textbf{A}$=\left\{\alpha>-1,\beta>\alpha,\beta>-\alpha\right\}$, it turns out that we can choose any
$$
0\leq\tau\leq1,
$$
and the matrix $\bm P$ in \eqref{Pst} will always be stochastic. If $\alpha\geq-1/2$ then this is possible for any $\gamma>-1,$ while for $-1<\alpha<-1/2$ we need to have $\gamma+\alpha+3/2>0$.
\item In the region \textbf{B}$=\left\{\beta>-1,\alpha>\beta\right\}$, it turns out that we can choose any
$$
0\leq\tau\leq C_{-1/2},
$$
and the matrix $\bm P$ in \eqref{Pst} will always be stochastic. If $\beta\geq-1/2$ then this is possible for any $\gamma>-1,$ while for $-1<\beta<-1/2$ we need to have $\gamma+\beta+3/2>0$.
\item In the region \textbf{C}$=\left\{\alpha>-1,\beta>\alpha,\beta<-\alpha\right\}$, it turns out that we can choose any
$$
0\leq\tau\leq\min\{C_{1/2},C_{\gamma+1}\},
$$
and the matrix $\bm P$ in \eqref{Pst} will always be stochastic. If $\alpha\geq-1/2$ then this is possible for any $\gamma>-1,$ while for $-1<\alpha<-1/2$ we need to have $\gamma+\alpha+3/2>0$. This is the only case where the upper bound may depend on $\gamma$. A straightforward computation shows that $C_{1/2}\leq C_{\gamma+1}$ if $\gamma\geq-1/2$.
\end{itemize}

Therefore, for all values of $\tau$ in the ranges described above for the regions \textbf{A}, \textbf{B} and \textbf{C}, we have a \emph{family} of discrete-time QBD processes $\{Z_t : t = 0, 1, \ldots\}$ with transition probability matrix $\bm P=(1-\tau)J_1+\tau J_2$. Thus the Karlin-McGregor representation formula (see formula (2.13) of \cite{FdI})  for the $(i,j)$ block entry of the matrix $\bm P$ is given by
\begin{equation}\label{KMC1}
\bm P_{i,j}^n=\left(\int_\Omega[(1-\tau)u+\tau v]^n\mathbb{Q}_i(u,v)\mathbb{Q}_j^T(u,v)W_{\alpha,\beta,\gamma}(u,v)dudv\right)\Pi_j,
\end{equation}
where $\mathbb{Q}_n,n\geq0,$ are the vector polynomials satisfying \eqref{TTRRQ} and $W_{\alpha,\beta,\gamma}$ is the normalized weight function \eqref{wKpJ2}. The matrices $\Pi_j, j\geq0,$ are the inverses of the norms of the corresponding vector polynomials $\mathbb{Q}_j, j\geq0$. Using \cite[Lemma 2.1]{FdI} we have one way of giving an explicit expression of these norms (another way could be using \cite[Section 6]{Sp76}). Indeed, it is possible to see that a generalized inverse $G_n=(G_{n,1},G_{n,2})$ of $C_n^T=(C_{n,1},C_{n,2})^T$ is given by
$$
G_n=\left[\begin{array}{ccccc|ccc}
1/c_{n,0} & 0&\cdots&0&0&0&\cdots&0  \\
0 & 1/c_{n,1}&\cdots&0&0&0&\cdots&0 \\
\vdots & \vdots&\ddots&\vdots&\vdots&\vdots&&\vdots\\
0 &0&\cdots&0&1/c_{n,n-1}&0&\cdots&0\\
0 &0&\cdots&-\frac{c_{n,n-2}^{(3)}}{c_{n,n}^{(1)}c_{n,n-2}}&-\frac{c_{n,n-1}^{(2)}}{c_{n,n}^{(1)}c_{n,n-1}}&0&\cdots&1/c_{n,n}^{(1)}
\end{array}\right].
$$
Since the representation of $\Pi_n$ is independent of the choice of the generalized inverse $G_n$ (see \cite[Lemma 2.1]{FdI}) we have, after some straightforward computations, that $\Pi_n$ is a diagonal matrix of the form $\Pi_n=\mbox{diag}\left(\Pi_{n,0},\cdots,\Pi_{n,n}\right)$ where
\begin{align*}
\Pi_{n,k}&=\frac{(\alpha+1)_k(\alpha+\gamma+3/2)_n(\alpha+\beta+\gamma+5/2)_{n-1}(\alpha+\beta+2\gamma+3)_{n+k-1}(2\gamma+2)_{n-k-1}(\alpha+\beta+2)_{k-1}}{(\beta+1)_k(\beta+\gamma+3/2)_n(\alpha+\beta+2)_{n+k}(\gamma+3/2)_n}\\
&\qquad\times\frac{(2n+\alpha+\beta+2\gamma+2)(n+k+\alpha+\beta+\gamma+3/2)(2n-2k+2\gamma+1)(2k+\alpha+\beta+1)}{k!(n-k)!}.
\end{align*}
Another way of writing these norms in terms of the norms of the Jacobi polynomials \eqref{normJ} is
\begin{align*}
\Pi_{n,k}&=\|Q_{n+\gamma+1/2}^{(\beta,\alpha)}\|^{-2}\|Q_{k}^{(\beta,\alpha)}\|^{-2}\|Q_{\gamma+1/2}^{(\beta,\alpha)}\|^{2}\\
&\qquad\times\frac{(\alpha+\beta+2\gamma+2)_{n+k}(2\gamma+2)_{n-k-1}(n+k+\alpha+\beta+\gamma+3/2)(2n-2k+2\gamma+1)}{(\alpha+\beta+2)_{n+k}(n-k)!(\alpha+\beta+\gamma+3/2)}.
\end{align*}
In particular we have that the family of polynomials $\mathbb{Q}_n,n\geq0,$ is \emph{mutually orthogonal}. Therefore, another way to write the Karlin-McGregor \eqref{KMC1} formula is entry by entry
\begin{align*}
\left(\bm P_{i,j}^n\right)_{i',j'}&=\mathbb{P}\left[Z_t=(j,j')\,|\,Z_0=(i,i')\right)]=\frac{\Pi_{j,j'}}{C\sigma_{i,i'}\sigma_{j,j'}}\sum_{k=0}^n\binom{n}{k}(1-\tau)^{n-k}\tau^k\\
\times&\left(\int_\Omega u^{n-k}v^{k}P_{i,i'}^{\alpha,\beta,\gamma}(u,v)P_{j,j'}^{\alpha,\beta,\gamma}(u,v)(1-2u+v)^\alpha(2u+v-1)^\beta(2u^2-2u-v+1)^\gamma dudv\right).
\end{align*}
According to \cite[Theorem 2.5]{FdI} we can construct an invariant measure $\bm\pi$ for the QBD process which is given by
\begin{align*}
\bm\pi=&\left(\Pi_0;\left(\Pi_1\bm e_2\right)^T;\left(\Pi_3\bm e_3\right)^T;\cdots\right)\\
=&\left(1; \frac{(2\alpha+2\gamma+3)(\alpha+\beta+2\gamma+4)(2\alpha+2\beta+2\gamma+5)}{(\alpha+\beta+2)(2\beta+2\gamma+3)},\right.\\
&\qquad\left. \frac{(\alpha+1)(2\alpha+2\gamma+3)(2\alpha+2\beta+2\gamma+7)(\alpha+\beta+2\gamma+3)_2}{(\beta+1)(2\beta+2\gamma+3)(2\gamma+3)(\alpha+\beta+2)};\cdots\right).
\end{align*}
Here $\bm e_N$ denotes the $N$-dimensional vector with all components equal to 1. Finally, it is also possible to study recurrence of the family of discrete-time QBD processes using (2.21) of \cite{FdI}. The process is recurrent if and only if
$$
\int_\Omega\frac{W_{\alpha,\beta,\gamma}(u,v)}{1-\tau v-(1-\tau)u}dudv=\infty.
$$
After some computations it turns out that, in the range of the values of $\tau$ for which $\bm P$ is stochastic, this integral is divergent, and therefore (null) recurrent, if and only if $-3/2<\alpha+\gamma\leq-1$. Otherwise the QBD process is transient. The QBD process can never be positive recurrent since the spectral measure is absolutely continuous and does not have any jumps (see the end of Section 2 of \cite{FdI} for more details).

\begin{remark}
Instead of \eqref{Pst}, we could have considered the situation where $\bm P=\tau_1J_1+\tau_2J_2$ and $\bm P\bm e=\bm 0$, in which case we would have had a continuous-time QBD process. Since $J_i\bm e=\bm e, i=1,2$ then we need $\tau_2=-\tau_1=-\tau$ and therefore $\bm P=\tau(J_2-J_1)$. All off-diagonal entries of $\bm P$ must be nonnegative while the entries of the main diagonal must be nonpositive. A closer look to these conditions entry by entry shows that it is never possible to have a continuous-time QBD process in this context.
\end{remark}

\begin{remark}
Going back to Remark \ref{rem32} we could have studied under what conditions we may provide a probabilistic interpretation of a linear combination of $J_1$ and $J_2$ of the form $\bm P=\tau_1J_1+\tau_2J_2$ for the case where we normalize the polynomials at the point $(0,1)$. For that there are at least two possibilities, either a continuous or a discrete-time QBD process. If we want to have a continuous-time QBD process then we need $\bm P\bm e = \bm0$ and nonnegative off-diagonal entries. But this is possible if and only if $\tau_2=0$ and $\tau_1>0$, i.e. a positive scalar multiple of $J_1$. If we want to have a discrete-time QBD process then we need $\bm P\bm e = \bm e$ and nonnegative (scalar) entries. This is possible if and only if $\tau_2=1$ and the parameter $\tau_1=\tau$ is chosen in such a way that all entries of $\bm P$ are nonnegative. The entries of $\bm P=\tau J_1+J_2$ are nonnegative if and only if
\begin{align*}
a_{n,k}^{(1)}>0&, \quad\tau a_{n,k}+a_{n,k}^{(2)}>0, \quad a_{n,k}^{(3)}>0,\\
\tau d_{n,k}+b_{n,k}^{(1)}\geq0&,  \quad\tau(b_{n,k}-1)+b_{n,k}^{(2)}\geq0,  \quad\tau e_{n,k}+b_{n,k}^{(3)}\geq0,\\
c_{n,k}^{(1)}>0&,  \quad\tau c_{n,k}+c_{n,k}^{(2)}>0,  \quad c_{n,k}^{(3)}>0.
\end{align*}
Now, from the definition (see \eqref{coeffpar0} and \eqref{coeffpar}), we have
$$
\frac{a_{n,k}^{(2)}}{a_{n,k}}=\frac{c_{n,k}^{(2)}}{c_{n,k}}=2(2b_k-1),\quad \frac{b_{n,k}^{(1)}}{d_{n,k}}=\frac{b_{n,k}^{(3)}}{e_{n,k}}=2(2b_{n+\gamma+1/2}-1).
$$
Therefore, as before, the lower bounds for $\tau$ (depending also on $\alpha,\beta,\gamma$) will probably depend on the behavior of the constant value $C_\gamma=2(1-2b_{\gamma+1/2})$. Additionally, the condition $\tau(b_{n,k}-1)+b_{n,k}^{(2)}\geq0$ is equivalent to $\tau\leq b_{n,k}^{(2)}/(1-b_{n,k})$, meaning that will also have upper bounds for $\tau$. We leave the details to the reader.
\end{remark}

\section{An urn model for the Jacobi-Koornwinder bivariate polynomials}\label{sec4}

In this section we will give an urn model associated with one of the QBD models introduced in the previous section. For simplicity, we will study the case of the discrete-time QBD process \eqref{Pst} with $\tau=0$ (therefore $\bm P=J_1$) and  $\beta=\alpha$. In this section we will assume that $\alpha$ and $\gamma$ are nonnegative integers. Consider $\{Z_t : t=0,1,\ldots\}$ the discrete-time QBD process on the state space $\{(n,k) : 0\leq k\leq n,n\in\mathbb{N}_0\}$ whose one-step transition probability matrix is given by the coefficients $A_{n,1}, B_{n,1}, C_{n,1}$ in \eqref{ABC1}-\eqref{coeffpar0} (see also \eqref{coefTTRR}). At every time step $t = 0, 1, 2,\ldots$ the state $(n, k)$ will represent the number of $n$ blue balls inside the $k$-th urn A$_k, k = 0,1,\ldots,n$. Observe that the number of urns available goes with the number of blue balls at every time step. All the urns we use sit in a bath consisting of an infinite number of blue and red balls.

\begin{figure}[h]
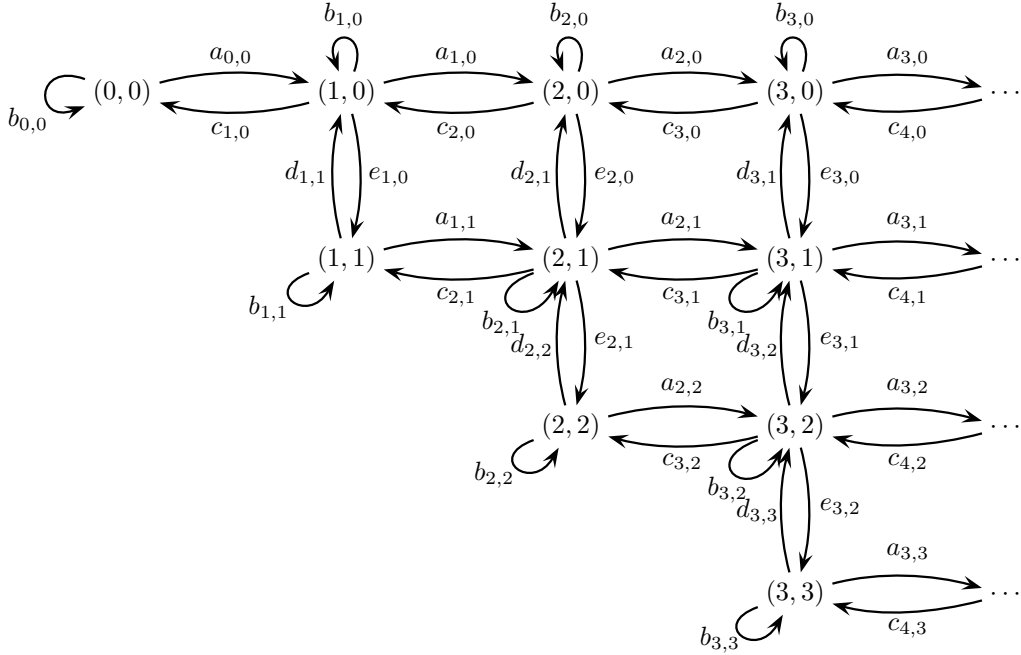

\vspace{.75cm}\begin{center}
$$\begin{psmatrix}[rowsep=1.8cm,colsep=2.2cm]
  \rnode{0}{(0,0)}& \rnode{1}{(1,0)} & \rnode{3}{(2,0)}& \rnode{6}{(3,0)}& \rnode{10}{\Huge{\cdots}} \\
    & \rnode{2}{(1,1)} & \rnode{4}{(2,1)}& \rnode{7}{(3,1)}& \rnode{11}{\Huge{\cdots}} \\
     & & \rnode{5}{(2,2)}& \rnode{8}{(3,2)}& \rnode{12}{\Huge{\cdots}} \\
    & & & \rnode{9}{(3,3)}& \rnode{13}{\Huge{\cdots}}
\psset{nodesep=3pt,arcangle=15,labelsep=2ex,linewidth=0.3mm,arrows=->,arrowsize=1mm
3} \nccurve[angleA=160,angleB=200,ncurv=4]{0}{0}
\nccurve[angleA=70,angleB=110,ncurv=6]{1}{1}
\nccurve[angleA=70,angleB=110,ncurv=6]{3}{3}
\nccurve[angleA=70,angleB=110,ncurv=6]{6}{6}
\nccurve[angleA=200,angleB=240,ncurv=4]{2}{2}
\nccurve[angleA=200,angleB=240,ncurv=4]{5}{5}
\nccurve[angleA=200,angleB=240,ncurv=4]{9}{9}
\nccurve[angleA=200,angleB=240,ncurv=5]{4}{4}
\nccurve[angleA=200,angleB=240,ncurv=5]{8}{8}
\nccurve[angleA=200,angleB=240,ncurv=5]{7}{7}
 \ncarc{0}{1}\ncarc{1}{0} \ncarc{1}{2} \ncarc{2}{1}
\ncarc{3}{4} \ncarc{4}{3} \ncarc{4}{5}\ncarc{5}{4} \ncarc{7}{6} \ncarc{6}{7}
\ncarc{3}{1}\ncarc{1}{3}\ncarc{3}{6}\ncarc{6}{3}\ncarc{6}{10}\ncarc{10}{6}\ncarc{2}{4}\ncarc{4}{2}\ncarc{4}{7}\ncarc{7}{4}\ncarc{7}{11}\ncarc{11}{7}\ncarc{5}{8}\ncarc{8}{5}\ncarc{8}{12}\ncarc{12}{8}\ncarc{9}{13}\ncarc{13}{9}\ncarc{7}{8}\ncarc{8}{7}\ncarc{9}{8}\ncarc{8}{9}
\uput[u](-13.2,5.9){b_{0,0}}\uput[u](-10.0,3.4){b_{1,1}}\uput[u](-7.0,1.2){b_{2,2}}
\uput[u](-4.0,-1.0){b_{3,3}}
\uput[u](-9.0,7.3){b_{1,0}}\uput[u](-6.0,7.3){b_{2,0}}\uput[u](-3.0,7.3){b_{3,0}}
\uput[u](-6.9,3.2){b_{2,1}}\uput[u](-3.9,3.2){b_{3,1}}\uput[u](-3.9,1.0){b_{3,2}}
\uput[u](-10.5,6.8){a_{0,0}}\uput[u](-7.5,6.8){a_{1,0}}\uput[u](-4.5,6.8){a_{2,0}}\uput[u](-1.5,6.8){a_{3,0}}
\uput[u](-7.5,4.6){a_{1,1}}\uput[u](-4.5,4.6){a_{2,1}}\uput[u](-1.5,4.6){a_{3,1}}
\uput[u](-4.5,2.4){a_{2,2}}\uput[u](-1.5,2.4){a_{3,2}}\uput[u](-1.5,0.2){a_{3,3}}
\uput[u](-10.5,5.8){c_{1,0}}\uput[u](-7.5,5.8){c_{2,0}}\uput[u](-4.5,5.8){c_{3,0}}\uput[u](-1.5,5.8){c_{4,0}}
\uput[u](-7.5,3.6){c_{2,1}}\uput[u](-4.5,3.6){c_{3,1}}\uput[u](-1.5,3.6){c_{4,1}}
\uput[u](-4.5,1.4){c_{3,2}}\uput[u](-1.5,1.4){c_{4,2}}\uput[u](-1.5,-0.8){c_{4,3}}
\uput[u](-8.4,5.2){e_{1,0}}\uput[u](-5.4,5.2){e_{2,0}}\uput[u](-2.4,5.2){e_{3,0}}\uput[u](-5.4,3.0){e_{2,1}}\uput[u](-2.4,3.0){e_{3,1}}\uput[u](-2.4,0.8){e_{3,2}}
\uput[u](-9.5,5.2){d_{1,1}}\uput[u](-6.5,5.2){d_{2,1}}\uput[u](-3.5,5.2){d_{3,1}}\uput[u](-6.5,2.9){d_{2,2}}\uput[u](-3.5,2.9){d_{3,2}}\uput[u](-3.5,0.7){d_{3,3}}
\end{psmatrix}
$$
\end{center}
\vspace{.5cm}
\caption{Diagram of all possible transitions of the discrete-time QBD process corresponding with $J_1$ for the bivariate Jacobi-Koornwinder polynomials.}
\label{fig2}
\end{figure}

Since $\beta=\alpha$ the coefficients in \eqref{coeffpar0} are simplified and given explicitly by
\begin{equation}\label{coeffpar00}
\begin{split}
a_{n,k}&=\frac{(2n+4\alpha+2\gamma+3)(n-k+2\gamma+1)(n+k+2\alpha+2\gamma+ 2)}{4(n+\alpha+\gamma+1)(2n-2k+2\gamma+1)(2n+2k+4\alpha+2\gamma+3)},\quad k=0,1,\ldots,n, \\
c_{n,k}&=\frac{(2n+2\gamma+1)(n-k)(n+k+2\alpha+1)}{4(n+\alpha+\gamma+1)(2n-2k+2\gamma+1)(2n+2k+4\alpha+2\gamma+3)},\quad k=0,1,\ldots,n-1,
\\
e_{n,k}&=\frac{(k+2\alpha+1)(n-k)(n+k+2\alpha+2\gamma+2)}{(2k+2\alpha+1)(2n-2k+2\gamma+1)(2n+2k+4\alpha+2\gamma+3)},\quad k=0,1,\ldots,n-1,\\
d_{n,k}&=\frac{k(n-k+2\gamma+1)(n+k+ 2\alpha+1)}{(2k+2\alpha+1)(2n-2k+2\gamma+1)(2n+2k+4\alpha+2\gamma+3)},\quad k=1,2,\ldots,n-1,\\
b_{n,k}&=\frac{1}{2},\quad k=0,1,\ldots,n.
\end{split}
\end{equation}
In Figure \ref{fig2} we can see a diagram of all possible transitions of this discrete-time QBD process.

At time $t=0$ the initial state is $Z_0=(n,k)$. The urn model will be divided in two steps. First, we consider two auxiliary urns U$_1$ and U$_2$. In urn U$_1$ we put $n+k+2\alpha+2\gamma+2$ blue balls and $n+k+2\alpha+1$ red balls from the bath, and in urn U$_2$ we put $n-k+2\gamma+1$ blue balls and $n-k$ red balls also from the bath. Then we draw independently one ball from urn U$_1$ and urn U$_2$ at random with the uniform distribution. We have four possibilities:

\begin{enumerate}
\item Both balls from U$_1$ and U$_2$ are blue with probability
$$
\frac{n+k+2\alpha+2\gamma+2}{2n+2k+4\alpha+2\gamma+3}\times\frac{n-k+2\gamma+1}{2n-2k+2\gamma+1}.
$$
Observe that this number is included in the coefficient $a_{n,k}$ in \eqref{coeffpar00}. Then we remove all the balls in urn A$_k$ and put $2n+4\alpha+2\gamma+3$ blue balls and $2n+2\gamma+1$ red balls in urn A$_k$. Draw one ball from A$_k$ at random with the uniform distribution. If we get a blue ball then we remove all balls in urn A$_k$ and add $n+1$ blue balls to the urn A$_k$ and start over. Therefore we have
$$
\mathbb{P}\left[Z_{1}=(n+1,k)\, |\, Z_0=(n,k)\right]=a_{n,k}.
$$
\item Both balls from U$_1$ and U$_2$ are red with probability
$$
\frac{n+k+2\alpha+1}{2n+2k+4\alpha+2\gamma+3}\times\frac{n-k}{2n-2k+2\gamma+1}.
$$
Observe that this number is included in the coefficient $c_{n,k}$ in \eqref{coeffpar00}. Then we remove all the balls in urn A$_k$ and put $2n+4\alpha+2\gamma+3$ blue balls and $2n+2\gamma+1$ red balls in urn A$_k$. Draw one ball from A$_k$ at random with the uniform distribution. If we get a red ball then we remove all balls in urn A$_k$ and add $n-1$ blue balls to the urn A$_k$ and start over. Therefore we have
$$
\mathbb{P}\left[Z_{1}=(n-1,k)\, |\, Z_0=(n,k)\right]=c_{n,k}.
$$
\item The ball from U$_1$ is blue and the ball from U$_2$ is red with probability
$$
\frac{n+k+2\alpha+2\gamma+2}{2n+2k+4\alpha+2\gamma+3}\times\frac{n-k}{2n-2k+2\gamma+1}.
$$
Observe that this number is included in the coefficient $e_{n,k}$ in \eqref{coeffpar00}. Then we remove all the balls in urn A$_k$ and put $k+2\alpha+1$ blue balls and $k$ red balls in urn A$_k$. Draw one ball from A$_k$ at random with the uniform distribution. If we get a blue ball then we remove all balls in urn A$_k$ and add $n$ blue balls to the urn A$_{k+1}$ and start over. Therefore we have
$$
\mathbb{P}\left[Z_{1}=(n,k+1)\, |\, Z_0=(n,k)\right]=e_{n,k}.
$$
\item The ball from U$_1$ is red and the ball from U$_2$ is blue with probability
$$
\frac{n+k+2\alpha+1}{2n+2k+4\alpha+2\gamma+3}\times\frac{n-k+2\gamma+1}{2n-2k+2\gamma+1}.
$$
Observe that this number is included in the coefficient $d_{n,k}$ in \eqref{coeffpar00}. Then we remove all the balls in urn A$_k$ and put $k+2\alpha+1$ blue balls and $k$ red balls in urn A$_k$. Draw one ball from A$_k$ at random with the uniform distribution. If we get a red ball then we remove all balls in urn A$_k$ and add $n$ blue balls to the urn A$_{k-1}$ and start over. Therefore we have
$$
\mathbb{P}\left[Z_{1}=(n,k-1)\, |\, Z_0=(n,k)\right]=d_{n,k}.
$$
\end{enumerate}
In each of the previous four possibilities there is a complementary probability. In cases (1) and (3) we may have a red ball in the second step while in cases (2) and (4) we may have a blue ball in the second step. In all these four possibilities we remove all balls in urn A$_k$ and add $n$ blue balls to the urn A$_{k}$ and start over. The addition of these four probabilities gives $1/2$. Therefore we have
$$
\mathbb{P}\left[Z_{1}=(n,k)\, |\, Z_0=(n,k)\right]=b_{n,k}=\frac{1}{2}.
$$

\begin{remark}
If $\beta\neq\alpha$ the probabilities in \eqref{coeffpar0} will have an extra factor, so we will have to add an extra step to the previous urn model. However, the urn model is not as clear as the previous one.
\end{remark}

\begin{remark}
It would be possible to consider an urn model taking $\tau=1$ in \eqref{Pst} (therefore $\bm P=J_2$). But in this case the coefficients $A_{n,2}, B_{n,2}, C_{n,2}$ in \eqref{ABC2}-\eqref{coeffpar} are way more complicated than the case we studied here. The diagram of all possible transitions will look like Figure 3 of \cite{FdI}. In \cite{FdI} an urn model was proposed for the orthogonal polynomials on the triangle as a consequence of finding a simple stochastic LU factorization of $\bm P$. Although it may be possible to find a LU factorization of $\bm P$ in this situation, each of the factors are not as simple as the original one, so this method is no longer convenient to find an urn model.
\end{remark}

\end{document}